# ZERO TEMPERATURE LIMIT FOR INTERACTING BROWNIAN PARTICLES. II. COAGULATION IN ONE DIMENSION[1]

By Tadahisa Funaki

*University of Tokyo*

We study the zero temperature limit for interacting Brownian particles in one dimension with a pairwise potential which is of finite range and attains a unique minimum when the distance of two particles becomes $a > 0$. We say a chain is formed when the particles are arranged in an "almost equal" distance $a$. If a chain is formed at time 0, so is for positive time as the temperature of the system decreases to 0 and, under a suitable macroscopic space-time scaling, the center of mass of the chain performs the Brownian motion with the speed inversely proportional to the total mass. If there are two chains, they independently move until the time when they meet. Then, they immediately coalesce and continue the evolution as a single chain. This can be extended for finitely many chains.

**1. Introduction.** We consider a system of interacting Brownian particles in a real line $\mathbb{R}$. The positions of $N$ particles at time $t$ are denoted by $\mathbf{x}(t) = (x_i(t))_{i=1}^N \in \mathbb{R}^N$ and evolve according to the stochastic differential equation (SDE)

$$(1.1) \qquad dx_i(t) = -\frac{1}{2}\varepsilon^{-\alpha}\frac{\partial H}{\partial x_i}(\mathbf{x}(t))\,dt + dw_i(t), \qquad 1 \le i \le N.$$

Here $(w_i(t))_{i=1}^N$ is a family of independent one-dimensional standard Brownian motions. The parameter $\varepsilon > 0$, which is very small, represents the ratio of the microscopic spatial unit length to the macroscopic one and $\varepsilon^{-\alpha}$ with $\alpha > 0$ is the inverse temperature of the system, which is already rescaled in $\varepsilon$. The Hamiltonian $H(\mathbf{x})$ of the configuration $\mathbf{x} = (x_i)_{i=1}^N \in \mathbb{R}^N$ is defined as a sum of pairwise interactions between particles:

$$(1.2) \qquad H(\mathbf{x}) = \sum_{1 \le i < j \le N} U(x_i - x_j).$$

---

Received November 2002; revised April 2003.

[1]Supported in part by JSPS Grants (B)(1) 14340029 and 13874015.

The potential $U = U(|x|)$ is symmetric, smooth, of finite range and has a unique nondegenerate minimum at $|x| = a > 0$; see Assumptions I and II stated in Sections 2 and 3 for details. The configuration $\mathbf{x} = (x_i)_{i=1}^N$ is a microscopic object and its macroscopic correspondence is given by $(\varepsilon x_i)_{i=1}^N$ under the spatial scaling $x \mapsto \varepsilon x$.

The aim of this paper is to investigate the asymptotic behavior of the system as $\varepsilon \downarrow 0$ under a proper scaling both in particles' number $N$ and time $t$ besides the temperature of the system. The particles' number $N \equiv N(\varepsilon)$ will change in such a manner that

$$(1.3) \qquad \lim_{\varepsilon \downarrow 0} \varepsilon N(\varepsilon) = \rho$$

with some $\rho > 0$. We shall simply write $N \sim \rho \varepsilon^{-1}$ for (1.3). Then, the solution $\mathbf{x}(t) = (x_i(t))_{i=1}^{N(\varepsilon)}$ of the SDE (1.1) is rescaled in time as

$$(1.4) \qquad \mathbf{x}^{(\varepsilon)}(t) = \mathbf{x}(\varepsilon^{-3}t), \qquad t \geq 0.$$

As the temperature, given by $\varepsilon^\alpha$, tends to 0, the system of the particles is expected to be frozen and arranged in an almost equal distance $a$. This naturally leads us to the following notion: a configuration $\mathbf{x} = (x_i)_{i=1}^N$ arranged in increasing order is called a *chain* with fluctuation $c \geq 0$ if it satisfies $|x_{i+1} - x_i - a| \leq c$ for every $1 \leq i \leq N - 1$. When $\varepsilon \downarrow 0$, the fluctuation $c \equiv c(\varepsilon)$ of the chain is expected to be small. Macroscopically, a *rod*, which is an interval $[\varepsilon x_1, \varepsilon x_{N(\varepsilon)}]$ in $\mathbb{R}$, rather than a set of points $(\varepsilon x_i)_{i=1}^{N(\varepsilon)}$ is associated with the chain $\mathbf{x} = (x_i)_{i=1}^{N(\varepsilon)}$ with small fluctuation under the spatial scaling $x \mapsto \varepsilon x$. The constants $\rho$ and $\rho a$ represent the mass and the length of the associated rod, respectively.

The first paper [1] studied the behavior of the rescaled process $\mathbf{x}^{(\varepsilon)}(t)$ for a single crystal which is an extended notion of chain in higher dimensions, and the result can be reformulated as follows in one dimension. If $\mathbf{x}^{(\varepsilon)}(0)$ is a chain with particles' number $N \sim \rho \varepsilon^{-1}$ and fluctuation $\varepsilon^{\nu^*}$ with certain $\nu^* > 0$ at $t = 0$, then $\mathbf{x}^{(\varepsilon)}(t) = (x_i^{(\varepsilon)}(t))_{i=1}^N$ remains to be a chain with fluctuation $\varepsilon^\nu$ with slightly smaller $\nu$ than $\nu^*$ for $t > 0$ asymptotically with probability one as $\varepsilon \downarrow 0$. Moreover, the (macroscopic) center of mass of the associated rod defined by $\eta^{(\varepsilon)}(t) := \frac{\varepsilon}{N} \sum_{i=1}^N x_i^{(\varepsilon)}(t)$ behaves asymptotically as $\varepsilon \downarrow 0$ as $\eta(0) + w(t/\rho)$ if $\eta(0) = \lim_{\varepsilon \downarrow 0} \eta^{(\varepsilon)}(0)$ exists, where $w(t)$ is the one-dimensional standard Brownian motion. This means that the evolutional speed of the rod is proportional to the inverse of the macroscopic mass $\rho$, see Theorem 2.2 for details.

The main result of this paper is stated in Section 3. Assume that $\mathbf{x}^{(\varepsilon)}(0) = \mathbf{x}^{(\varepsilon,1)}(0) \cup \mathbf{x}^{(\varepsilon,2)}(0)$ consists of two chains $\mathbf{x}^{(\varepsilon,1)}(0)$ and $\mathbf{x}^{(\varepsilon,2)}(0)$ with $N_1 \sim \rho_1 \varepsilon^{-1}$ and $N_2 \sim \rho_2 \varepsilon^{-1}$ particles, respectively, where $\rho_1, \rho_2 > 0$. Then, these two chains evolve independently until they meet. Once they meet, they immediately coalesce and form a larger single chain with particles' number $N = N_1 + N_2$, see Theorem 3.1. Afterwards it evolves as a single chain, so the associated rod performs the Brownian motion with speed inversely proportional to $\rho_1 + \rho_2$; see Corollary 3.9. This result can be extended for finitely many chains; see Corollary 3.9.

As is explained in [1], one of the motivations of our study comes from the theory of interfaces and, in this respect, the rod we have introduced can be regarded as a kind of Wulff shape at temperature zero. The system of sticky Brownian motions was discussed by Smoluchowski; see [6]. Another model for coalescing rods in one dimension was studied by Mullins [7]. These models have, however, a slightly different character from ours, since the speeds of the particles or the rods do not change after coalescence in these models. Lang [5] investigated a system of ordinary differential equations (1.1) dropping Brownian motions with $\varepsilon^{-\alpha} = 1$ and $N = \infty$. Such system arises from the SDE in the zero temperature limit under a proper time change.

**2. Motion of a single chain.** This section summarizes the results for a single chain, which are deduced from Theorem 3.4 of Funaki [1] by restricting the system in one dimension.

2.1. *Hamiltonian.* The space $\mathbb{R}_*^N$ stands for the set of all $\mathbf{x} = (x_i)_{i=1}^N \in \mathbb{R}^N$ arranged in increasing order $x_1 \leq x_2 \leq \cdots \leq x_N$. The Hamiltonian $H(\mathbf{x})$ of $\mathbf{x}$ is introduced by the formula (1.2). The pair potential $U$ in (1.2) satisfies the following conditions:

ASSUMPTION I. (i) (symmetry) $U(x) = U(-x)$, $x \in \mathbb{R}$.

(ii) (smoothness, finite range) $U \in C_0^3(\mathbb{R})$.

(iii) There exists a unique $a > 0$ such that $U(a) = \min_{x \geq 0} U(x)$ and $\check{c} := U''(a) > 0$.

(iv) $b < 2a$, where $b := \inf\{x > 0; U(y) = 0$ for every $y > x\}$.

We denote by $\mathbf{z}$ the configuration such that $z_{i+1} - z_i = a, 1 \leq i \leq N - 1$. Note that $\mathbf{z}$ is a local minimum, which is sometimes called an instanton in physics, of the Hamiltonian $H$. By Assumption I(iv), each particle in the configuration $\mathbf{z}$ interacts only with neighboring particles. The (microscopic) center of mass of the configuration $\mathbf{x} \in \mathbb{R}_*^N$ is defined by

$$(2.1) \qquad \eta(\mathbf{x}) = \frac{1}{N} \sum_{i=1}^N x_i \in \mathbb{R}.$$

Let $\mathbf{z}^0 = (z_i^0)_{i=1}^N$ be the centered local minimum $\mathbf{z}$, that is, $\eta(\mathbf{z}^0) = 0$. Then, each configuration $\mathbf{x} \in \mathbb{R}_*^N$ can be decomposed as

$$(2.2) \qquad \mathbf{x} = \mathbf{z}^0 + \mathbf{h} + \eta(\mathbf{x})$$

with $\mathbf{h} \equiv \mathbf{h}(\mathbf{x}) = (h_i)_{i=1}^N \in \mathbb{R}^N$ satisfying $\sum_{i=1}^N h_i = 0$, where $\mathbf{z}^0 + \mathbf{h} + \eta :=
(z_i^0 + h_i + \eta)_{i=1}^N$ for $\eta \in \mathbb{R}$. For $\mathbf{h} = (h_i)_{i=1}^N \in \mathbb{R}^N$ satisfying $\sum_{i=1}^N h_i = 0$, we
introduce three norms $\|\nabla \mathbf{h}\|_2, \|\nabla \mathbf{h}\|_\infty$ and $\|\Delta \mathbf{h}\|_2$, respectively, by

$$\|\nabla \mathbf{h}\|_2^2 = \sum_{i=1}^{N-1} (h_{i+1} - h_i)^2, \qquad \|\nabla \mathbf{h}\|_\infty = \max_{1 \le i \le N-1} |h_{i+1} - h_i|,$$

$$\|\Delta \mathbf{h}\|_2^2 = \sum_{i=2}^{N-1} (2h_i - h_{i+1} - h_{i-1})^2 + (h_2 - h_1)^2 + (h_N - h_{N-1})^2.$$

LEMMA 2.1.  *For every* $\mathbf{h} \in \mathbb{R}^N$,

$$\tfrac{1}{2}\|\Delta \mathbf{h}\|_2 \le \|\nabla \mathbf{h}\|_2 \le N\|\Delta \mathbf{h}\|_2.$$

PROOF.   The first inequality is obvious from the definition of two norms.
To show the next, set $g_i = h_{i+1} - h_i, 1 \le i \le N-1$. Then, we have

$$\|\nabla \mathbf{h}\|_2^2 = \sum_{i=1}^{N-1} g_i^2 = \sum_{i=1}^{N-1} \left\{ g_1 + \sum_{j=2}^i (g_j - g_{j-1}) \right\}^2$$

$$\le \sum_{i=1}^{N-1} i \left\{ g_1^2 + \sum_{j=2}^i (g_j - g_{j-1})^2 \right\} \le \frac{N(N-1)}{2}\|\Delta \mathbf{h}\|_2^2.$$

This implies the second inequality.  □

Two quadratic forms $\mathcal{E}_1$ and $\mathcal{E}_2$ of $\mathbf{h}$ introduced in Lemmas 2.1 and 3.1
of [1] related to the Hamiltonian $H(\mathbf{x})$ have the forms $\mathcal{E}_1(\mathbf{h}) = \check{c}\|\nabla \mathbf{h}\|_2^2$ and
$\mathcal{E}_2(\mathbf{h}) = \check{c}^2\|\Delta \mathbf{h}\|_2^2$ in one dimension, respectively. Lemma 2.1 shows that the
constant $\lambda^{(2)}(\mathbf{z})$ arising in a bound between these two quadratic forms stated
in Lemma 3.2 of [1] can be taken as $\lambda^{(2)}(\mathbf{z}) = \check{c}N^{-2}$.

Let $\mathcal{M} \equiv \mathcal{M}^N = \{\mathbf{z}^0 + \eta; \eta \in \mathbb{R}\}$ be the set of local minima and let $\mathcal{M}^\nabla(c) \equiv
\mathcal{M}^{\nabla,N}(c) = \{\mathbf{x} \in \mathbb{R}_*^N; \|\nabla \mathbf{h}(\mathbf{x})\|_\infty \le c\}$ be the tubular neighborhood of $\mathcal{M}$ for
$c \in [0, b-a]$. The configuration $\mathbf{x} \in \mathcal{M}^\nabla(c)$ will be called a *chain* with parti-
cles' number $N$ and fluctuation $c$. Note that $\mathbf{x} \in \mathcal{M}^\nabla(c)$ means $|x_{i+1} - x_i -
a| \le c$ for every $1 \le i \le N-1$.

2.2. *Microscopic shape theorem and motion of the macroscopic center
of mass.*  We now discuss the scaling limit for the solution $\mathbf{x}(t)$ of the
SDE (1.1). The particles' number of the system is assumed to behave as
$N \equiv N(\varepsilon) \sim \rho\varepsilon^{-1}$ with $\rho > 0$. Let $\mathbf{x}^{(\varepsilon)}(t) = (x_i^{(\varepsilon)}(t))_{i=1}^{N(\varepsilon)} \in \mathbb{R}^{N(\varepsilon)}$ be the time
changed process of $\mathbf{x}(t)$ defined by (1.4). For $\nu > 0$, consider the stopping
time $\sigma \equiv \sigma^{(\varepsilon)}$ determined by

$$\sigma = \inf \{t \ge 0; \mathbf{x}^{(\varepsilon)}(t) \notin \mathcal{M}^{\nabla,N(\varepsilon)}(\varepsilon^\nu)\} \equiv \inf \{t \ge 0; \|\nabla \mathbf{h}(\mathbf{x}^{(\varepsilon)}(t))\|_\infty > \varepsilon^\nu\}.$$

THEOREM 2.2. (1) *Assume* $\nu > 2, \alpha > 2\nu + 3$ *and* $\mathbf{x}^{(\varepsilon)}(0) \in \mathcal{M}^{\nabla, N(\varepsilon)}(\bar{c}(\varepsilon))$ *with* $\bar{c}(\varepsilon) = o(\varepsilon^{\nu+1/2})$ *as* $\varepsilon \downarrow 0$. *Then, we have* $\lim_{\varepsilon \downarrow 0} P(\sigma^{(\varepsilon)} \geq t) = 1$ *for every* $t > 0$.

(2) *Let* $\eta^{(\varepsilon)}(t) := \varepsilon \eta(\mathbf{x}^{(\varepsilon)}(t))$ *be the (macroscopic) center of mass of the rod associated with the chain* $\mathbf{x}^{(\varepsilon)}(t)$. *Then,* $\eta^{(\varepsilon)}(t)$ *weakly converges to* $\eta(0) + w(t/\rho)$ *as* $\varepsilon \downarrow 0$ *in the space* $C([0, T], \mathbb{R})$ *for every* $T > 0$ *if* $\eta(0) = \lim_{\varepsilon \downarrow 0} \eta^{(\varepsilon)}(0)$ *exists, where* $w(t)$ *is the one-dimensional standard Brownian motion.*

PROOF. We apply Theorem 3.4 of [1] to show assertion (1). Note that $\lambda^{(1,\varepsilon)} \equiv \check{c}$ (see Remark 2.2 of [1]), $\lambda^{(2,\varepsilon)} = \check{c}N^{-2} \sim \check{c}\rho^{-2}\varepsilon^2$ (see the remark after Lemma 2.1) in one dimension and we take $c(\varepsilon) = \varepsilon^\nu$. Then, condition (3.6) in [1] is satisfied if $c(\varepsilon) \leq c\varepsilon^2$ (as $\varepsilon \downarrow 0$) for some small $c > 0$, which is valid when $\nu > 2$. (The condition $c(\varepsilon) \leq \bar{c}(\varepsilon)$ in (3.6) automatically holds for sufficiently small $\varepsilon > 0$, since we can take $\bar{c}(\varepsilon) = b - a$; see Remark 3.2(1) of [1] as well.) To see that condition (3.7) in [1] holds, note that

$$\|\nabla \mathbf{h}(\mathbf{x}^{(\varepsilon)}(0))\|_2^2 \leq N \|\nabla \mathbf{h}(\mathbf{x}^{(\varepsilon)}(0))\|_\infty^2 \leq C\varepsilon^{-1}\{o(\varepsilon^{\nu+1/2})\}^2 = o(\varepsilon^{2\nu})$$

and therefore

$$\{\lambda^{(1,\varepsilon)}c(\varepsilon)^2\}^{-p} E[\|\nabla \mathbf{h}(\mathbf{x}^{(\varepsilon)}(0))\|_2^{2p}] \leq C\varepsilon^{-2\nu p} \times o(\varepsilon^{2\nu p}) \to 0.$$

Condition (3.8) in [1] also holds since

$$\{\lambda^{(1,\varepsilon)}c(\varepsilon)^2\}^{-p}\varepsilon^{-\kappa}\beta(\varepsilon)^{-p+1}N(\varepsilon)^p\{\lambda^{(2,\varepsilon)}\}^{-p+1} \leq C\varepsilon^{(\alpha-2\nu-3)p-\alpha-1} \to 0$$

for large $p$; recall $\kappa = 3, \beta(\varepsilon) = \varepsilon^{-\alpha}$ and $\alpha - 2\nu - 3 > 0$ from our assumption. Assertion (2) is easy, since $\eta^{(\varepsilon)}(t) = \eta^{(\varepsilon)}(0) + \varepsilon N^{-1}\sum_{i=1}^N w_i(\varepsilon^{-3}t)$. The constants $C$ in the above estimates may change from line to line. $\square$

Theorem 2.2(1) asserts that asymptotically with probability one $\mathbf{x}^{(\varepsilon)}(t)$ remains to be a chain with fluctuation $\varepsilon^\nu$ if it is a chain with fluctuation $o(\varepsilon^{\nu+1/2})$ at $t = 0$. This characterizes the microscopic structure of the solutions of the SDE (1.1) which are scaled macroscopically in time. Theorem 2.2(2) determines the macroscopic evolution of the associated rod.

**3. Coagulation of two chains.** In this section, we assume that there are two chains $\mathbf{x}^{(1)} \equiv \mathbf{x}^{(\varepsilon,1)}$ and $\mathbf{x}^{(2)} \equiv \mathbf{x}^{(\varepsilon,2)}$ in $\mathbb{R}$ with particles' numbers $N_1 \equiv N_1(\varepsilon)$ and $N_2 \equiv N_2(\varepsilon)$, respectively. The chain $\mathbf{x}^{(1)}$ is located on the left side of $\mathbf{x}^{(2)}$. If the distance between the right most particle of $\mathbf{x}^{(1)}$ and the left most one of $\mathbf{x}^{(2)}$ (which will be called the distance of two chains) is larger than $b$, these two chains move independently. We shall show that, once the distance of two chains becomes $b$, these two chains coalesce immediately in macroscopic time scale and afterwards move as a single chain with particles' number $N \equiv N(\varepsilon) := N_1 + N_2$.

To be more precise, we assume the following conditions on $\mathbf{x}(0)$ [$= \mathbf{x}^{(\varepsilon)}(0)$, a sequence depending on $\varepsilon > 0$] and consider the solution $\mathbf{x}(t) = \mathbf{x}^{(1)}(t) \cup \mathbf{x}^{(2)}(t)$ of the SDE (1.1) starting at $\mathbf{x}(0)$, where $\mathbf{x}^{(1)}(t) = (x_i(t))_{i=1}^{N_1}$ and $\mathbf{x}^{(2)}(t) = (x_i(t))_{i=N_1+1}^{N}$.

CONDITION A.   (i) $\mathbf{x}(0) = \mathbf{x}^{(1)}(0) \cup \mathbf{x}^{(2)}(0)$ consists of two chains $\mathbf{x}^{(1)}(0) = (x_i(0))_{i=1}^{N_1}$ and $\mathbf{x}^{(2)}(0) = (x_i(0))_{i=N_1+1}^{N}$ with particles' numbers $N_1 \sim \rho_1 \varepsilon^{-1}, N_2 \sim \rho_2 \varepsilon^{-1}$ and fluctuation $\varepsilon^\mu, \mu > 0$, that is, $\mathbf{x}^{(\ell)}(0) \in \mathcal{M}^{\nabla, N_\ell}(\varepsilon^\mu), \ell = 1, 2$, where $\rho_1, \rho_2 > 0$.

(ii) The distance of these two chains is $b$, that is, $x_{N_1+1}(0) - x_{N_1}(0) = b$.

We need, in addition to Assumption I, the following rather technical assumptions on the shape of the potential $U$. We shall denote the connected component of the set $\{x > 0; U''(x) > 0\}$ containing $a$ by $D = (b_1, b_2)$.

ASSUMPTION II.   (i) $2U(b_2) > U(a)$ and $U(b_1) + U(b_2) > U(a)$.

(ii) $U'(x) \geq 0$ for every $x \geq b_2$ (and therefore for $x \geq a$).

(iii) $2b_3 > b$, where $b_3 \in (b_1, a)$ is determined by $U(b_2) + U(b_3) = U(a)$.

Note that, since $U(a) < U(a) - U(b_2) < U(b_1)$ from (i) and (ii), $b_3$ in (iii) exists uniquely. Assumptions (i) and (iii) mean that the well at $a$ is deep and located away from 0, respectively. An example of the potential $U$, which satisfies both Assumptions I and II, is given by $U(x) = \psi((|x| - a)^2 - 4)$, where we assume $a \geq 4$ and $\psi \in C^3(\mathbb{R})$ is a nondecreasing function such that $\psi(x) = x$ for $x \leq -1$ and $\psi(x) = 0$ for $x \geq 0$. Note that $U(b_1), U(b_2) > -1, b \leq a + 2$ and $b_3 > a - 1$ in this example.

The main result of this section is now formulated.

THEOREM 3.1.   Let $\mathbf{x}^{(\varepsilon)}(t) = \mathbf{x}(\varepsilon^{-3}t)$ be the time changed process of $\mathbf{x}(t)$ with initial data $\mathbf{x}(0)$ satisfying Condition A with $\mu > 1/2$. Assume that $\tilde{\nu} > 0$ is given and $\alpha$ satisfies $\alpha > 4 \vee (2\tilde{\nu} + 3)$. Then we have, for every $\delta > 0$,

$$\lim_{\varepsilon \downarrow 0} P(\mathbf{x}^{(\varepsilon)}(t) \in \mathcal{M}^{\nabla, N(\varepsilon)}(\varepsilon^{\tilde{\nu}}) \text{ for some } t \leq \varepsilon^{1-\delta}) = 1.$$

Theorem 3.1 combined with Theorem 2.2 establishes the asymptotic behavior of two chains located in a general position. Suppose that $\nu > 2, \tilde{\nu} > 5/2$ are given and $\alpha > (2\nu + 3) \vee (2\tilde{\nu} + 3)$, and that the initial data $\mathbf{x}(0) = \mathbf{x}^{(1)}(0) \cup \mathbf{x}^{(2)}(0)$ of the SDE (1.1) satisfies only Condition A(i) for some $\mu > \nu + 1/2$. Then, by Theorem 2.2, two chains $\mathbf{x}^{(\varepsilon, \ell)}(t), \ell = 1, 2$, scaled macroscopically in time both stay in $\mathcal{M}^{\nabla, N_\ell}(\varepsilon^\nu)$ until the time when the distance of these two chains becomes $b$. However, at the time when the distance of two chains becomes $b$, Condition A(ii) is also satisfied and therefore we can apply Theorem 3.1 to see that, within the time $\varepsilon^{1-\delta}$, a single chain $\mathbf{x}^{(\varepsilon)}(t) \in \mathcal{M}^{\nabla,N}(\varepsilon^{\tilde{\nu}})$ is formed from two chains $\mathbf{x}^{(\varepsilon,\ell)}(t), \ell = 1, 2$. Afterward, applying Theorem 2.2 again, the single chain $\mathbf{x}^{(\varepsilon)}(t)$ moves staying in $\bigcap_{\delta > 0} \mathcal{M}^{\nabla,N}(\varepsilon^{\tilde{\nu}-1/2-\delta})$. All these statements hold asymptotically with probability one as $\varepsilon \downarrow 0$. The coagulation of several chains will be discussed in Corollary 3.9.

The first step for the proof of Theorem 3.1 is to show that, asymptotically with probability one as $\varepsilon \downarrow 0$, the distances of all neighboring particles of $\mathbf{x}(t)$ belong to the convex region $D''$ of the potential $U$ [which is slightly smaller than $D$; see (3.22)] at certain time $t$ smaller than $\varepsilon^{-2-\delta}$ for arbitrary $\delta > 0$, see Proposition 3.5 and the remark after it. In the proof, we always consider $\mathbf{x}(t)$ without introducing the scaling in time. The second step is to prove that, once the distances of all neighboring particles belong to $D''$, the two chains coagulate and form a single chain within the time $\varepsilon^{\alpha-2-\delta}$ for every $\delta > 0$; see Proposition 3.8.

Let $\mathbf{z}^{(1,2)} = \mathbf{z}^{(1)} \cup \mathbf{z}^{(2)}, \mathbf{z}^{(1)} = (z_i)_{i=1}^{N_1}, \mathbf{z}^{(2)} = (z_i)_{i=N_1+1}^{N}$ be the configuration satisfying $z_{i+1} - z_i = a, i \neq N_1$ and $z_{N_1+1} - z_{N_1} = b$. Note that $\mathbf{z}^{(1,2)}$ is a saddle point of $H(\mathbf{x})$. Condition A means that $\mathbf{x}(0)$ is in a neighborhood of $\mathbf{z}^{(1,2)}$.

From Assumptions II(i) and II(iii), there exists $b_4 \in (a, b_2)$ such that $U(b_4) = U(b_3)$. Then, choose an interval $D' = (b_3', b_4')$ and $b_2'$ in such a manner that $b_1 < b_3' < b_3, b_4 < b_4' < b_2 < b_2, 2b_3' > b, U(b_3') = U(b_4')$ and $\bar{\delta} := U(b_2') + U(b_3') - U(a) > 0$. This is possible by taking $b_2', b_3'$ slightly smaller than $b_2, b_3$ and $b_4'$ slightly larger than $b_4$, respectively. We introduce four stopping times:

$$\tau_1 = \inf\{t \geq 0; x_{N_1+1}(t) - x_{N_1}(t) \leq b_2'\},$$

$$\tau_2 = \inf\{t \geq 0; x_{i+1}(t) - x_i(t) \notin D' \text{ for some } i \neq N_1\},$$

$$\tau_3 = \inf\{t \geq 0; \tilde{H}(\mathbf{x}(t)) \geq \delta_1\},$$

$$\tau_4 = \inf\left\{t \geq 0; \eta(\mathbf{x}^{(2)}(t)) - \eta(\mathbf{x}^{(1)}(t)) \leq \frac{a}{2}N - N^\kappa\right\},$$

where $\tilde{H}(\mathbf{x}) := H(\mathbf{x}) - H(\mathbf{z}^{(1,2)}), 0 < \delta_1 < \bar{\delta}$ and $0 < \kappa < 1$; $\kappa$ will be chosen later in the proof of Proposition 3.5. The functions $\eta(\mathbf{x}^{(1)}) = \frac{1}{N_1}\sum_{i=1}^{N_1} x_i$ and $\eta(\mathbf{x}^{(2)}) = \frac{1}{N_2}\sum_{i=N_1+1}^{N} x_i$, defined by (2.1) with particles' number $N$ replaced by $N_1$ and $N_2$ in each chain, represent the microscopic centers of mass of $\mathbf{x}^{(1)}$ and $\mathbf{x}^{(2)}$, respectively. We first discuss with $\tau_1, \tau_2, \tau_3$ in Lemmas 3.2, 3.3 and $\tau_4$ will be treated in Lemma 3.4, later; see Remark 3.1 for the meaning of $\tau_4$. These three lemmas are prepared for the proof of Proposition 3.5.

LEMMA 3.2. (i) *For every* $t \leq \tau_1 \wedge \tau_2$, $\min_{1 \leq i \leq N-2}\{x_{i+2}(t) - x_i(t)\} \geq b$. *In particular, in the configuration* $\mathbf{x} := \mathbf{x}(t)$, *only neighboring particles interact and we have*

$$(3.1) \quad \tilde{H}(\mathbf{x}) = \sum_{1 \leq i \leq N-1, i \neq N_1}\{U(x_{i+1} - x_i) - U(a)\} + U(x_{N_1+1} - x_{N_1}).$$

(ii) $\tau_1 \wedge \tau_3 < \tau_2$.

PROOF. For $t \leq \tau_1 \wedge \tau_2$, we have $x_{N_1+1}(t) - x_{N_1}(t) \geq b_2'$ and $b_4' \leq x_{i+1}(t) - x_i(t) \leq b_4'$ for every $i \neq N_1$. This implies $x_{i+2}(t) - x_i(t) \geq 2b_3'(\geq b)$ for all $1 \leq i \leq N-2$, and therefore (i) is shown. To prove (ii), assume $\tau_2 \leq \tau_1 \wedge \tau_3$ and set $\mathbf{x} \equiv (x_i)_{i=1}^{N} := \mathbf{x}(\tau_2)$. Then, there exists $i_0 \ (\neq N_1)$ such that $x_{i_0+1} - x_{i_0} = b_3'$ (or $b_4'$), and $x_{N_1+1} - x_{N_1} \geq b_2'$ because $\tau_2 \leq \tau_1$. Moreover, since $\tau_2 \leq \tau_1$, we can apply (1) at $t = \tau_2$ and see that $\tilde{H}(\mathbf{x})$ has the form (3.1). Therefore,

$$\tilde{H}(\mathbf{x}) \geq \{U(x_{i_0+1} - x_{i_0}) - U(a)\} + U(x_{N_1+1} - x_{N_1})$$
$$\geq \{U(b_3') - U(a)\} + U(b_2') = \bar{\delta} > \delta_1.$$

This contradicts $\tau_2 \leq \tau_3$ and the proof of (ii) is complete. □

LEMMA 3.3. *Assume* $\alpha > 4$ *and* $\mu > 1/2$ *in Condition* A *on* $\mathbf{x}(0)$. *Then, for every* $\gamma > 0$,

$$\lim_{\varepsilon \downarrow 0} P(\tau_1 \wedge \varepsilon^{-\gamma} \leq \tau_3) = 1.$$

PROOF. *Step* 1. Take $\delta_0 \in (0, \delta_1)$ and fix it. In this step, we prove that there exist $c_1, \varepsilon_0 > 0$ such that

$$(3.2) \quad \sum_{j=1}^{N}\left(\frac{\partial H}{\partial x_j}\right)^2(\mathbf{x}(t)) \geq \frac{c_1 \delta_0}{N^3}$$

if $t \leq \tau_1 \wedge \tau_2$, $\tilde{H}(\mathbf{x}(t)) \geq \delta_0$ and $\varepsilon \in (0, \varepsilon_0)$. Indeed, since (3.1) holds for $\mathbf{x} := \mathbf{x}(t)$ and since $x_{N_1+1} - x_{N_1} \geq b_2'$ implies $U(x_{N_1+1} - x_{N_1}) \leq 0$ [by Assumption II(ii)], we see from $\tilde{H}(\mathbf{x}) \geq \delta_0$ that $U(x_{i_0+1} - x_{i_0}) - U(a) \geq \delta_0/N$ holds for some $i_0 \neq N_1$. However, for such $i_0$, $|(x_{i_0+1} - x_{i_0}) - a| \geq c_2\sqrt{\delta_0/N}$ [by Assumption I(iii)] and accordingly $|U'(x_{i_0+1} - x_{i_0})| \geq c_3\sqrt{\delta_0/N}$ (by noting $x_{i_0+1} - x_{i_0} \in D'$) for certain $c_2, c_3 > 0$.

Now let us assume that (3.2) does not hold. Then, we have

$$(3.3) \quad \left|\frac{\partial H}{\partial x_j}(\mathbf{x})\right| = |U'(x_j - x_{j-1}) - U'(x_{j+1} - x_j)| \leq \sqrt{\frac{c_1 \delta_0}{N^3}}$$

for every $1 \leq j \leq N$; we regard as $U'(x_1 - x_0) = U'(x_{N+1} - x_N) = 0$. First consider the case where $i_0 \leq N_1 - 1$ and $U'(x_{i_0+1} - x_{i_0}) \geq c_3\sqrt{\delta_0/N}$. [The case where $i_0 \geq N_1 + 1$ or $U'(x_{i_0+1} - x_{i_0}) \leq -c_3\sqrt{\delta_0/N}$ can be similarly treated.] Then, using (3.3) with $j = i_0$, we have

$$U'(x_{i_0} - x_{i_0-1}) \geq c_3\sqrt{\delta_0/N} - \sqrt{c_1\delta_0/N^3}.$$

Continuing this procedure of estimates $i_0 - 1$ times, we finally arrive at

$$(3.4) \quad U'(x_2 - x_1) \geq c_3\sqrt{\delta_0/N} - (i_0 - 1)\sqrt{c_1\delta_0/N^3} \geq (c_3 - \sqrt{c_1})\sqrt{\delta_0/N}.$$

But, if one takes $c_1 > 0$ such that $c_3 > \sqrt{c_1}$, (3.4) contradicts (3.3) with $j = 1$ for sufficiently small $\varepsilon > 0$. Therefore (3.2) is shown.

*Step* 2. Simple application of Itô's formula for the solution $\mathbf{x}(t)$ of (1.1) shows

$$d\tilde{H}(\mathbf{x}(t)) = dm(t) + \{-\varepsilon^{-\alpha}b^{(1)}(\mathbf{x}(t)) + b^{(2)}(\mathbf{x}(t))\}\, dt,$$

where $b^{(1)}(\mathbf{x}) = \frac{1}{2}\sum_{j=1}^{N}(\partial H/\partial x_j)^2(\mathbf{x})$, $b^{(2)}(\mathbf{x}) = \frac{1}{2}\sum_{j=1}^{N}\partial^2 H/\partial x_j^2(\mathbf{x})$ and $m(t)$ is a martingale defined by

$$m(t) = \sum_{j=1}^{N}\int_0^t \frac{\partial H}{\partial x_j}(\mathbf{x}(s))\, dw_j(s).$$

However, in Step 1, we have seen $b^{(1)}(\mathbf{x}(t)) \geq c_1\delta_0/(2N^3) \geq c_4\varepsilon^3$ for some $c_4 > 0$ if $t \leq \tau_1 \wedge \tau_2$ and $\tilde{H}(\mathbf{x}(t)) \geq \delta_0$. Moreover, since $\partial^2 H/\partial x_j^2$ are bounded, $|b^{(2)}(\mathbf{x})| \leq c_5\varepsilon^{-1}$. Therefore, recalling $\alpha > 4$, we obtain

$$(3.5) \quad d\tilde{H}(\mathbf{x}(t)) \leq dm(t) - c_6\varepsilon^{-\alpha+3}\, dt,$$

for $t \leq \tau_1 \wedge \tau_2$ satisfying $\tilde{H}(\mathbf{x}(t)) \geq \delta_0$; or, more precisely saying, $\tilde{H}(\mathbf{x}(t)) - m(t)$ is differentiable in $t$ and $d\{\tilde{H}(\mathbf{x}(t)) - m(t)\}/dt \leq -c_6\varepsilon^{-\alpha+3}$ for such $t$. Since $\partial H/\partial x_j$ are bounded, the derivative of the quadratic variational process of $m(t)$ is dominated by

$$(3.6) \quad \frac{d}{dt}\langle m\rangle_t = \sum_{j=1}^{N}\left(\frac{\partial H}{\partial x_j}\right)^2(\mathbf{x}(t)) \leq c_7\varepsilon^{-1}.$$

*Step* 3. Introduce a time changed process $y_t$ of $\tilde{H}(\mathbf{x}(t))$ as $y_t := \tilde{H}(\mathbf{x}(\langle m\rangle_t^{-1}))$, where $\langle m\rangle_t^{-1}$ denotes the inverse function of $\langle m\rangle_t$. Then, from (3.5) and (3.6), we have

$$(3.7) \quad dy_t \leq dB_t - c\varepsilon^{-\alpha+4}\, dt$$

if $t \leq \langle m\rangle_{\tau_1 \wedge \tau_2}$ and $y_t \in [\delta_0, \infty)$, where $c := c_6 c_7^{-1}$. Note that $B_t := m(\langle m\rangle_t^{-1})$ is a Brownian motion and $y_0 = \tilde{H}(\mathbf{x}(0)) \leq c_8\varepsilon^{2\mu-1}$ from Condition A on $\mathbf{x}(0)$ and the bound $H(\mathbf{x}^{(\ell)}(0)) - H(\mathbf{z}^{(\ell)}) \leq C\mathcal{E}_1(\mathbf{h}(\mathbf{x}^{(\ell)}(0)))$, $\ell = 1, 2$, shown in (3.3) of [1], where $H(\mathbf{x}^{(\ell)})$ denotes the Hamiltonian of the system with $N_\ell$ particles.

Choose $\delta_2 \in (\delta_0, \delta_1)$ and take a smooth function $f : (\delta_0, \infty) \to [0, \infty)$, which satisfies $f(x) = 0$ for every $x \geq \delta_2$ and $s(\delta_0+) \equiv \lim_{x \downarrow \delta_0} s(x) = -\infty$, where $s(x) \equiv s^{(\varepsilon)}(x)$ is a function defined by

$$s(x) = \int_{\delta_2}^{x} \exp\left\{-2 \int_{\delta_2}^{y} (-c\varepsilon^{-\alpha+4} + f(z)) \, dz\right\} dy, \qquad x > \delta_0.$$

In fact, such function $f$ can be taken, since, if $f(x)$ behaves as $f(x) \sim C(x - \delta_0)^{-\lambda}$ as $x \downarrow \delta_0$ with $\lambda > 1$ and $C > 0$, then $s^{(\varepsilon)}(\delta_0+) = -\infty$ for each $\varepsilon > 0$. We consider the SDE for $z_t \equiv z_t^{(\varepsilon)}$:

$$(3.8) \qquad dz_t = dB_t - c\varepsilon^{-\alpha+4} \, dt + f(z_t) \, dt, \qquad z_0 = \delta_2.$$

The function $s(x)$ is the so-called natural scale (or canonical scale) for the diffusion process $z_t$; see [3] or [4], page 339. Since $s(\delta_0+) = -\infty$, it holds that

$$(3.9) \qquad z_t > \delta_0, \qquad t \geq 0, \text{ a.s.}$$

Moreover, for every sufficiently small $\varepsilon > 0$,

$$(3.10) \qquad y_t \leq z_t$$

holds for every $0 \leq t \leq \langle m \rangle_{\tau_1 \wedge \tau_2}$. Indeed, since $\mu > 1/2$, $y_0 \leq c_8 \varepsilon^{2\mu-1} \leq z_0$ (as $\varepsilon \downarrow 0$) and therefore (3.10) is true at $t = 0$. If $y_t \leq \delta_0$, (3.10) automatically holds because $z_t > \delta_0$. Once $y_t$ moves into the interval $[\delta_0, \infty)$, one can apply the comparison theorem (see, e.g., [2]) between two processes $y_t$ and $z_t$ recalling (3.7), (3.8), and $f \geq 0$, and (3.10) is shown for all $t \leq \langle m \rangle_{\tau_1 \wedge \tau_2}$.

We now consider a stopping time $\sigma_1 \equiv \sigma_1^{(\varepsilon)} = \inf\{t \geq 0; z_t = \delta_1\}$ for the solution of the SDE (3.8). Then, (3.10) implies $\sigma_1 \leq \langle m \rangle_{\tau_3}$ if $\tau_3 \leq \tau_1 \wedge \tau_2$, which shows

$$\{\tau_3 < \tau_1 \wedge \tau_2 \wedge \varepsilon^{-\gamma}\} \subset \{\sigma_1 \leq \langle m \rangle_{\varepsilon^{-\gamma}}\} \subset \{\sigma_1 \leq c_7 \varepsilon^{-1-\gamma}\}.$$

The second inclusion follows from (3.6). Hence, if one can show

$$(3.11) \qquad \lim_{\varepsilon \downarrow 0} P(\sigma_1 \leq c_7 \varepsilon^{-1-\gamma}) = 0,$$

then we have $\lim_{\varepsilon \downarrow 0} P(\tau_3 < \tau_1 \wedge \tau_2 \wedge \varepsilon^{-\gamma}) = 0$. However, since

$$\{\tau_3 < \tau_1 \wedge \varepsilon^{-\gamma}\} \subset \{\tau_3 < \tau_1 \wedge \tau_2 \wedge \varepsilon^{-\gamma}\} \cup \{\tau_2 \leq \tau_1 \wedge \tau_3\},$$

Lemma 3.2(2) concludes the lemma.

*Step* 4. Only the proof of (3.11) is left. The argument of this step is rather standard. We introduce another stopping time $\sigma_2 \equiv \sigma_2^{(\varepsilon)} = \inf\{t \geq 0; z_t = \delta_3\}$ by choosing $\delta_3 \in (\delta_2, \delta_1)$. Then, one can find $\bar{\lambda} > 0$ such that

$$(3.12) \qquad P(\sigma_2 > 1) \geq \bar{\lambda}$$

for every $0 < \varepsilon < 1$. In fact, consider an SDE

$$d\tilde{z}_t = dB_t + f(\tilde{z}_t)\,dt, \qquad \tilde{z}_0 = z_0 = \delta_2.$$

Then, $z_t \leq \tilde{z}_t$ holds for every $t \geq 0$, which implies $\tilde{\sigma}_2 \leq \sigma_2^{(\varepsilon)}$ for $\tilde{\sigma}_2 = \inf\{t \geq 0; \tilde{z}_t = \delta_3\}$. We may therefore take $\bar{\lambda} := P(\tilde{\sigma}_2 > 1) > 0$.

Let $\{\sigma^{(k)} \equiv \sigma^{(k,\varepsilon)}\}_{k=0,1,2,\ldots}$ and $K \equiv K^{(\varepsilon)}$ be a sequence of stopping times and a random variable inductively defined by $\sigma^{(0)} = 0$ and for $k = 1, 2, \ldots$,

$$\sigma^{(2k-1)} = \inf\{t > \sigma^{(2k-2)}; z_t = \delta_3\},$$

$$\sigma^{(2k)} = \inf\{t > \sigma^{(2k-1)}; z_t \notin (\delta_2, \delta_1)\},$$

$$K = \inf\{k \geq 1; z_{\sigma^{(2k)}} = \delta_1\},$$

respectively. Then, $\{\sigma^{(k)}\}_{k=0,1,2,\ldots}$ and $K$ have the following four properties: (i) $\{\sigma^{(2k-1)} - \sigma^{(2k-2)}\}_{k=1,2,\ldots}$ is an independent system, (ii) the law of $\sigma^{(2k-1)} - \sigma^{(2k-2)}$ is identical to that of $\sigma_2$ for each $k = 1, 2, \ldots$, (iii) $\sigma_1 \geq \sum_{k=1}^{K}(\sigma^{(2k-1)} - \sigma^{(2k-2)})$ and (iv) $K - 1$ has the geometric distribution: $P(K - 1 = n) = pq^n, n = 0, 1, 2, \ldots$ with $p \equiv p^{(\varepsilon)} := P(K = 1)$ and $q = 1 - p$. Indeed, (i) is a consequence of the strong Markov property of $z_t$, while (ii) and (iii) are obvious. To see (iv), one may note that $\{\bar{z}_k := z_{\sigma^{(2k)}}\}_{k=0,1,2,\ldots}$ forms a two state Markov chain on the set $\{\delta_2, \delta_1\}$ with the transition probability $P(\bar{z}_1 = \delta_1|\bar{z}_0 = \delta_2) = p$ and $P(\bar{z}_1 = \delta_2|\bar{z}_0 = \delta_2) = q$. Furthermore, for every sufficiently small $\varepsilon > 0$,

$$(3.13) \qquad p \leq \exp\{-\bar{c}\varepsilon^{-\alpha+4}\},$$

where $\bar{c} = c(\delta_1 - \delta_3)$. In fact, this is shown noting that $p = \{s(\delta_3) - s(\delta_2)\}/\{s(\delta_1) - s(\delta_2)\}$ and the natural scale $s(x)$ is given by

$$s(x) = \frac{1}{2c}\varepsilon^{\alpha-4}\{\exp\{2c\varepsilon^{-\alpha+4}(x - \delta_2)\} - 1\}$$

for $x \geq \delta_2$; recall that $f(x) = 0$ if $x \geq \delta_2$.

With the choice of $K_0 \equiv K_0^{(\varepsilon)} := \varepsilon\exp\{\bar{c}\varepsilon^{-\alpha+4}\}$ ($\leq \varepsilon/p$), (3.13) and property (iv) show that

$$P(K \leq K_0) = 1 - (1 - p)^{K_0} \to 0, \qquad \varepsilon \downarrow 0.$$

Therefore, from property (iii), the proof of (3.11) is complete if one can prove

$$(3.14) \qquad \lim_{\varepsilon\downarrow 0} P\left(\sum_{k=1}^{K_0}(\sigma^{(2k-1)} - \sigma^{(2k-2)}) \leq c_7\varepsilon^{-1-\gamma}\right) = 0$$

for every $\gamma > 0$. Set $X_k = \mathbb{1}_{\{\sigma^{(2k-1)} - \sigma^{(2k-2)} > 1\}}$ for $k = 1, 2, \ldots$. Then, $\{X_k\}_{k=1,2,\ldots}$ is a sequence of independent and identically distributed random variables such that $\lambda \equiv \lambda^{(\varepsilon)} := P(X_1 = 1) \geq \bar{\lambda} > 0$ from (3.12) and the property (ii). Since

$$Y := \sum_{k=1}^{K_0} X_k \leq \sum_{k=1}^{K_0} (\sigma^{(2k-1)} - \sigma^{(2k-2)}),$$

(3.14) follows from $\lim_{\varepsilon \downarrow 0} P(Y \leq c_7 \varepsilon^{-1-\gamma}) = 0$. But, this is easy from $E[Y] = \lambda K_0, E[(Y - E[Y])^2] = \lambda(1-\lambda)K_0$ by applying Chebyshev's inequality.   □

LEMMA 3.4.   *Assume $\mu + \kappa > 1$ for the constants $\mu, \kappa > 0$ appearing in Condition A on $\mathbf{x}(0)$ and in the definition of $\tau_4$, respectively. Then, for every $\delta > 0$,*

$$\lim_{\varepsilon \downarrow 0} P(\tau_1 \wedge \tau_4 \leq \varepsilon^{-(1+2\kappa+\delta)}) = 1.$$

PROOF.   From the SDE (1.1), we have

$$d\eta(\mathbf{x}^{(1)}(t)) = -\frac{1}{2N_1} \varepsilon^{-\alpha} U'(x_{N_1} - x_{N_1+1})\, dt + \frac{1}{N_1} \sum_{i=1}^{N_1} dw_i(t)$$

$$\geq \frac{1}{N_1} \sum_{i=1}^{N_1} dw_i(t), \qquad t \leq \tau_1,$$

and a similar bound on $d\eta(\mathbf{x}^{(2)}(t))$ from above for $t \leq \tau_1$; recall that $U'(x) \geq 0$ for $x \geq a$ [by Assumption II(ii)] and the symmetry of $U$ [by Assumption I(i)]. Hence,

$$\eta(\mathbf{x}^{(2)}(t)) - \eta(\mathbf{x}^{(1)}(t)) \leq \eta(\mathbf{x}^{(2)}(0)) - \eta(\mathbf{x}^{(1)}(0)) + (N_1^{-1} + N_2^{-1})^{1/2} w(t)$$

for $t \leq \tau_1$, where

$$w(t) := (N_1^{-1} + N_2^{-1})^{-1/2} \left( \frac{1}{N_2} \sum_{i=N_1+1}^{N} w_i(t) - \frac{1}{N_1} \sum_{i=1}^{N_1} w_i(t) \right)$$

is a Brownian motion. Let us introduce a stopping time:

$$\sigma_3 = \inf \{ t \geq 0;\ (N_1^{-1} + N_2^{-1})^{1/2} w(t) \leq -(b-a) - \varepsilon^\mu N - N^\kappa \}.$$

Then, decomposing $\mathbf{x}^{(1)}$ and $\mathbf{x}^{(2)}$ as in (2.2), respectively, we have

$$(3.15)\qquad \begin{aligned} \eta(\mathbf{x}^{(1)}) &= \frac{a}{2}(1 - N_1) + \frac{1}{N_1} \sum_{i=1}^{N_1} (h_i - h_{N_1}) + x_{N_1}, \\ \eta(\mathbf{x}^{(2)}) &= \frac{a}{2}(N_2 - 1) + \frac{1}{N_2} \sum_{i=N_1+1}^{N} (h_i - h_{N_1+1}) + x_{N_1+1}, \end{aligned}$$

since $\sum_{i=1}^{N_1} h_i = \sum_{i=N_1+1}^{N} h_i = 0$ and $z_M^0 = -z_1^0 = a(M-1)/2$ for the centered local minimum $\mathbf{z}^0 = (z_i^0)_{i=1}^M$ with particles' number $M$ (we take $M = N_1, N_2$). Therefore, Condition A on $\mathbf{x}(0) = \mathbf{x}^{(1)}(0) \cup \mathbf{x}^{(2)}(0)$ implies

$$\left| \{\eta(\mathbf{x}^{(2)}(0)) - \eta(\mathbf{x}^{(1)}(0))\} - \frac{a}{2}N - (b-a) \right|$$

$$\leq \varepsilon^\mu \left( \frac{1}{N_1} \sum_{i=1}^{N_1} (N_1 - i) + \frac{1}{N_2} \sum_{i=N_1+1}^{N} (i - N_1 - 1) \right) \leq \varepsilon^\mu N,$$

from which we see $\tau_4 \leq \sigma_3$ if $\sigma_3 \leq \tau_1$. Accordingly, we have $\{\tau_1 \wedge \tau_4 > \varepsilon^{-(1+2\kappa+\delta)}\} \subset \{\sigma_3 > \varepsilon^{-(1+2\kappa+\delta)}\}$. However, since $\mu + \kappa > 1$ implies $\varepsilon^\mu N \ll N^\kappa$ as $\varepsilon \downarrow 0$, we see $(N_1^{-1} + N_2^{-1})^{-1/2}(b - a + \varepsilon^\mu N + N^\kappa) \leq c\varepsilon^{-(1/2+\kappa)}$ for some $c > 0$ so that $\sigma_3 \leq \tilde{\sigma}_3 := \inf\{t \geq 0; w(t) \leq -c\varepsilon^{-(1/2+\kappa)}\}$. The scaling invariance of the Brownian motion shows $\tilde{\sigma}_3 = \varepsilon^{-(1+2\kappa)}\bar{\sigma}_3$ in law, where $\bar{\sigma}_3 := \inf\{t \geq 0; w(t) \leq -c\}$. Therefore, we get

$$P(\tau_1 \wedge \tau_4 > \varepsilon^{-(1+2\kappa+\delta)}) \leq P(\tilde{\sigma}_3 > \varepsilon^{-(1+2\kappa+\delta)}) = P(\bar{\sigma}_3 > \varepsilon^{-\delta}) \to 0$$

as $\varepsilon \downarrow 0$, which completes the proof. □

REMARK 3.1. If $\mathbf{x} = \mathbf{z}^{(1)} \cup \mathbf{z}^{(2)} = (z_i)_{i=1}^{N_1} \cup (z_i)_{i=N_1+1}^{N}$ satisfies $z_{i+1} - z_i = a$ for all $i \neq N_1$, then (3.15) taking $h_i = 0$ for all $i$ shows $\eta(\mathbf{z}^{(2)}) - \eta(\mathbf{z}^{(1)}) = \frac{a}{2}N - a + z_{N_1+1} - z_{N_1}$. In particular, if the distance of $\mathbf{z}^{(1)}$ and $\mathbf{z}^{(2)}$ is $a$ (i.e., $\mathbf{x} \in \mathcal{M}^N$), then $\eta(\mathbf{z}^{(2)}) - \eta(\mathbf{z}^{(1)}) = \frac{a}{2}N$ [cf. this with $\eta(\mathbf{z}^{(2)}) - \eta(\mathbf{z}^{(1)}) = \frac{a}{2}N + (b-a)$ for $\mathbf{x} = \mathbf{z}^{(1,2)}$]. This may explain the meaning of the stopping time $\tau_4$. The randomness coming from the Brownian motions $(w_i(t))_{i=1}^N$ helps to make the distance between two chains shorter. Such effect was measured by the difference of the centers of mass of two chains.

The results obtained in Lemmas 3.2, 3.3 and 3.4 are summarized in the following proposition.

PROPOSITION 3.5. *Assume $\alpha > 4$ and $\mu > 1/2$. Then, for every $\delta > 0$,*

$$\lim_{\varepsilon \downarrow 0} P(\tau_1 \leq \varepsilon^{-2-\delta}, \tau_1 < \tau_2) = 1.$$

PROOF. Taking $\frac{1}{2} < \kappa < \frac{1+\delta}{2} \wedge 1$ and denoting $\delta' := 1 - 2\kappa + \delta > 0$ by $\delta$ again, we may prove that

(3.16)     $$\lim_{\varepsilon \downarrow 0} P(\tau_1 \leq \varepsilon^{-(1+2\kappa+\delta)}, \tau_1 < \tau_2) = 1, \qquad \delta > 0.$$

*Step* 1. We first note that (3.16) can be deduced from

(3.17)     $$\lim_{\varepsilon \downarrow 0} P(\tau_1 \leq \varepsilon^{-(1+2\kappa+\delta)}) = 1.$$

In fact, this is seen from Lemma 3.2(2) and Lemma 3.3 with $\gamma = 1 + 2\kappa + \delta$, since

$$\{\tau_1 \wedge \tau_3 < \tau_2\} \cap \{\tau_1 \wedge \varepsilon^{-(1+2\kappa+\delta)} \leq \tau_3\} \cap \{\tau_1 \leq \varepsilon^{-(1+2\kappa+\delta)}\}$$
$$\subset \{\tau_1 \leq \varepsilon^{-(1+2\kappa+\delta)}, \tau_1 < \tau_2\}.$$

We now give the proof of (3.17). Lemmas 3.3 and 3.4 show $\lim_{\varepsilon \downarrow 0} P(A^{(\varepsilon)}) = 1$ for $A^{(\varepsilon)} = \{\tau_1 \wedge \varepsilon^{-(1+2\kappa+\delta)} \leq \tau_3, \tau_1 \wedge \tau_4 \leq \varepsilon^{-(1+2\kappa+\delta)}\}$. Assume $\tau_4 \leq \tau_3$ and $\tau_1 > \varepsilon^{-(1+2\kappa+\delta)}$ on the event $A^{(\varepsilon)}$, and set $\mathbf{x} := \mathbf{x}(\tau_4) = \mathbf{x}^{(1)} \cup \mathbf{x}^{(2)} = (x_i)_{i=1}^N$. Then, we have $\tilde{H}(\mathbf{x}) \leq \delta_1$ and $x_{N_1+1} - x_{N_1} \geq b_2'$ since $\tau_4 < \tau_1$. Accordingly, noting $\tau_4 \leq \tau_1 \wedge \tau_3 \leq \tau_1 \wedge \tau_2$ by Lemma 3.2(2), we see that formula (3.1) holds for $\tilde{H}(\mathbf{x})$ and

$$(3.18) \qquad \sum_{1 \leq i \leq N-1, i \neq N_1} \{U(x_{i+1} - x_i) - U(a)\} = \tilde{H}(\mathbf{x}) - U(x_{N_1+1} - x_{N_1})$$
$$\leq \bar{\delta} - U(b_2') = U(b_3') - U(a).$$

On the other hand, $\mathbf{x}$ satisfies

$$(3.19) \qquad \eta(\mathbf{x}^{(2)}) - \eta(\mathbf{x}^{(1)}) = \frac{a}{2}N - N^{\kappa}.$$

We shall prove in Step 2 that (3.18) and (3.19) are incompatible. Once this is proved, we have $\tau_4 > \tau_3$ or $\tau_1 \leq \varepsilon^{-(1+2\kappa+\delta)}$ on $A^{(\varepsilon)}$ and this shows (3.17).

*Step* 2. Set $\bar{C} := U(b_3') - U(a) > 0$. Then, (3.18) implies $U(x_{i+1} - x_i) - U(a) \leq \bar{C}$ for every $i \neq N_1$ and, in particular, $x_{i+1} - x_i \in D$. However, $U$ is dominated from below by a quadratic function on $D$, that is, there exists $c_- > 0$ such that

$$c_-(g-a)^2 \leq U(g) - U(a) \qquad \text{if } g \in D.$$

Therefore, (3.18) shows that for $g_i := x_{i+1} - x_i, 1 \leq i \leq N-1$,

$$(3.20) \qquad \sum_{1 \leq i \leq N, i \neq N_1} (g_i - a)^2 \leq \frac{\bar{C}}{c_-}.$$

Next, noting that $x_j = x_{N_1} - \sum_{i=j}^{N_1-1} g_i$ for $1 \leq j \leq N_1$ and $x_j = x_{N_1+1} + \sum_{i=N_1+1}^{j-1} g_i$ for $N_1 + 1 \leq j \leq N$, we rewrite the difference of the centers of mass of two chains in terms of $\mathbf{g} = (g_i)_{i \neq N_1}$:

$$\eta(\mathbf{x}^{(2)}) - \eta(\mathbf{x}^{(1)}) = (x_{N_1+1} - x_{N_1}) + \frac{1}{N_1} \sum_{i=1}^{N_1-1} i g_i + \frac{1}{N_2} \sum_{i=N_1+1}^{N-1} (N-i) g_i.$$

Hence, recalling that $x_{N_1+1} - x_{N_1} \geq b_2'$, (3.19) implies

$$(3.21) \qquad F(\mathbf{g}) := \frac{1}{N_1} \sum_{i=1}^{N_1-1} i(g_i - a) + \frac{1}{N_2} \sum_{i=N_1+1}^{N-1} (N-i)(g_i - a)$$
$$\leq -N^{\kappa} - b_2' + a.$$

However, using Schwarz's inequality and (3.20), $|F(\mathbf{g})|$ is dominated by

$$
\begin{aligned}
|F(\mathbf{g})| \leq{}& \frac{1}{N_1}\left(\sum_{i=1}^{N_1-1} i^2\right)^{1/2}\left(\sum_{i=1}^{N_1-1} (g_i-a)^2\right)^{1/2} \\
&+ \frac{1}{N_2}\left(\sum_{i=N_1+1}^{N-1} (N-i)^2\right)^{1/2}\left(\sum_{i=N_1+1}^{N-1} (g_i-a)^2\right)^{1/2} \\
\leq{}& cN^{1/2}
\end{aligned}
$$

for some $c > 0$, which contradicts (3.21) since $\kappa > 1/2$. Therefore, (3.18) and (3.19) are incompatible. $\square$

If $\tau_1 < \tau_2$, the solution $\mathbf{x} := \mathbf{x}(\tau_1) = (x_i)_{i=1}^N$ of the SDE (1.1) at time $\tau_1$ satisfies $x_{N_1+1} - x_{N_1} = b_2'$ and $x_{i+1} - x_i \in D'$ for every $i \neq N_1$. In particular, it holds that

$$
(3.22) \qquad x_{i+1} - x_i \in D'' \qquad \text{for every } 1 \leq i \leq N-1,
$$

where $D'' := (b_3', b_2'] \Subset D = (b_1, b_2)$. Note that $c_* := \inf_{x \in D''} U''(x) > 0$.

We now move to the second stage. We begin with the investigation of the classical flow determined by the SDE (1.1) dropping the noise terms. Let $\bar{\mathbf{x}}(t) = (\bar{x}_i(t))_{i=1}^N$ be the solution of the ordinary differential equation (ODE)

$$
(3.23) \qquad \frac{d\bar{x}_i}{dt} = -\frac{1}{2}\varepsilon^{-\alpha}\frac{\partial H}{\partial x_i}(\bar{\mathbf{x}}), \qquad 1 \leq i \leq N,
$$

with an initial data $\bar{\mathbf{x}}(0) = \mathbf{x}$ satisfying the condition (3.22), and set

$$
g_i(t) = \bar{x}_{i+1}(t) - \bar{x}_i(t), \qquad 1 \leq i \leq N-1.
$$

Then, as long as $\min_{1 \leq i \leq N-2}\{\bar{x}_{i+2}(t) - \bar{x}_i(t)\} \geq b$, $\mathbf{g}(t) = (g_i(t))_{i=1}^{N-1}$ satisfies the ODE

$$
(3.24) \qquad \frac{dg_i}{dt} = \frac{1}{2}\varepsilon^{-\alpha}\{U'(g_{i+1}) + U'(g_{i-1}) - 2U'(g_i)\}, \qquad 1 \leq i \leq N-1,
$$

where $g_0(t) = g_N(t) := a$ in the right-hand side. The first assertion in the next lemma is the maximum principle, while the second is an energy inequality for the ODE (3.24). The convexity of $U$ on $D''$ is essential.

LEMMA 3.6. *Assume that $g_i(0) \in D''$ for all $1 \leq i \leq N-1$. Then, for every $t > 0$, we have*

$$
(3.25) \qquad g_i(t) \in D'' \qquad \text{for all } 1 \leq i \leq N-1,
$$

*and*

$$
(3.26) \qquad \sum_{i=1}^{N-1}(g_i(t)-a)^2 \leq \exp\{-c_*\varepsilon^{-\alpha}N^{-2}t\}\sum_{i=1}^{N-1}(g_i(0)-a)^2.
$$

Proof.  Assume that (3.25) holds at some $t \geq 0$. If $g_{i_0}(t) = \max_{0 \leq i \leq N} g_i(t)$ for such $t$ with some $1 \leq i_0 \leq N - 1$, then, since $U'$ is increasing on $D''$, the ODE (3.24) gives $dg_{i_0}(t)/dt \leq 0$ so that $g_{i_0}(t)$ is nonincreasing. Therefore, $\max_{0 \leq i \leq N} g_i(t)$ is also nonincreasing in $t$ [remembering the boundary conditions $g_0(t) = g_N(t) = a$]. Similarly, if $g_{i_0}(t) = \min_{0 \leq i \leq N} g_i(t)$ for some $1 \leq i_0 \leq N - 1$, then $g_{i_0}(t)$ and accordingly $\min_{0 \leq i \leq N} g_i(t)$ are nondecreasing. This shows that $g_i(t)$ can not go outside of $D''$ for all $1 \leq i \leq N - 1$. Thus assertion (3.25) is shown. To prove (3.26), we see from the ODE (3.24) that

$$
\begin{aligned}
\frac{d}{dt} \sum_{i=1}^{N-1} &(g_i(t) - a)^2 \\
&= -\varepsilon^{-\alpha} \sum_{i=0}^{N-1} (g_{i+1}(t) - g_i(t))\{U'(g_{i+1}(t)) - U'(g_i(t))\} \\
&\leq -c_* \varepsilon^{-\alpha} \sum_{i=0}^{N-1} (g_{i+1}(t) - g_i(t))^2 \\
&\leq -c_* \varepsilon^{-\alpha} N^{-2} \sum_{i=1}^{N-1} (g_i(t) - a)^2.
\end{aligned}
$$

(3.27)

The second line is from $U'' \geq c_*$ on $D''$, while the third line is by the Poincaré inequality: $\sum_{i=1}^{N-1} \bar{g}_i^2 \leq N^2 \sum_{i=0}^{N-1} (\bar{g}_{i+1} - \bar{g}_i)^2$ if $\bar{g}_0 = 0$, applied for $\bar{g}_i = g_i - a$. The bound (3.26) follows from (3.27). The Poincaré inequality is immediate from Schwarz's inequality as we saw in the proof of Lemma 2.1.  □

We shall next prove that, asymptotically with probability one, the solution $\mathbf{x}(t)$ of the SDE (1.1) moves along with the solution $\bar{\mathbf{x}}(t)$ of the ODE (3.23). This implies, with the help of Lemma 3.6, that $\mathbf{x}(t)$ goes into a neighborhood of a single chain; see Proposition 3.8.

Assume that $\mathbf{x}(t)$ and $\bar{\mathbf{x}}(t)$ have a common initial data $\mathbf{x} = \mathbf{x}(0) = \bar{\mathbf{x}}(0)$ satisfying condition (3.22) and introduce a stopping time:

$$
\tau_5 = \inf \left\{ t \geq 0; \max_{1 \leq i \leq N} |x_i(t) - \bar{x}_i(t)| \geq \varepsilon^\theta \right\}, \qquad \theta > 0.
$$

Lemma 3.7.  For every $\delta > 0$, we have

$$
\lim_{\varepsilon \downarrow 0} P(\tau_5 \geq \varepsilon^{2\theta + 1 + \delta}) = 1.
$$

Proof.  Applying Itô's formula for $I(t) := \sum_{i=1}^{N} (x_i(t) - \bar{x}_i(t))^2$, we have

$$
I(t) = I(0) + m(t) - \varepsilon^{-\alpha} \int_0^t b(\mathbf{x}(s), \bar{\mathbf{x}}(s)) \, ds + Nt,
$$

where

$$b(\mathbf{x}, \bar{\mathbf{x}}) = \sum_{i=1}^{N} (x_i - \bar{x}_i) \left\{ \frac{\partial H}{\partial x_i}(\mathbf{x}) - \frac{\partial H}{\partial x_i}(\bar{\mathbf{x}}) \right\},$$

$$m(t) = 2 \sum_{i=1}^{N} \int_0^t (x_i(s) - \bar{x}_i(s)) \, dw_i(s).$$

Denote the $2\varepsilon^\theta$-neighborhood of $D''$ by $D''_\varepsilon := (b'_3 - 2\varepsilon^\theta, b'_2 + 2\varepsilon^\theta)$. Then, since $\inf_{x \in D''_\varepsilon} U''(x) \geq 0$ and $2(b'_3 - 2\varepsilon^\theta) \geq b$ (for sufficiently small $\varepsilon > 0$), we have

$$b(\mathbf{x}, \bar{\mathbf{x}}) = \sum_{i=1}^{N-1} \left\{ (x_{i+1} - \bar{x}_{i+1}) - (x_i - \bar{x}_i) \right\} \left\{ U'(x_{i+1} - x_i) - U'(\bar{x}_{i+1} - \bar{x}_i) \right\} \geq 0,$$

if $x_{i+1} - x_i, \bar{x}_{i+1} - \bar{x}_i \in D''_\varepsilon$ for all $1 \leq i \leq N - 1$. Noting that (3.25) implies $\bar{x}_{i+1}(t) - \bar{x}_i(t) \in D''$ for every $t \geq 0$ and $1 \leq i \leq N - 1$, we see that $x_{i+1}(t) - x_i(t) \in D''_\varepsilon$ for every $t \leq \tau_5$ and $1 \leq i \leq N - 1$. Therefore, recalling $I(0) = 0$, we obtain $I(t) \leq m(t) + Nt$ for every $t \leq \tau_5$, and accordingly $E[I(\tau_5 \wedge t)] \leq Nt$ for all $t \geq 0$. Since $I(\tau_5) \geq \varepsilon^{2\theta}$ and $N \leq C\varepsilon^{-1}$, we have, for $\lambda = 2\theta + 1 + \delta$,

$$P(\tau_5 \leq \varepsilon^\lambda) \leq \varepsilon^{-2\theta} E[I(\tau_5 \wedge \varepsilon^\lambda)] \leq \varepsilon^{-2\theta} C\varepsilon^{-1} \varepsilon^\lambda = C\varepsilon^\delta \to 0, \qquad \varepsilon \downarrow 0. \quad \square$$

Lemmas 3.6 and 3.7 can be summarized into the following proposition for the stopping time $\tau \equiv \tau^{(\varepsilon)}$ defined by

$$(3.28) \qquad \tau = \inf \left\{ t > 0; \, \mathbf{x}(t) \in \mathcal{M}^{\nabla, N}(\varepsilon^{\tilde{\nu}}) \right\}, \qquad \tilde{\nu} > 0.$$

PROPOSITION 3.8. *Assume $\alpha > 2\tilde{\nu} + 3$ and $\mathbf{x}(0)$ satisfies condition* (3.22). *Then, for every $\delta > 0$,*

$$(3.29) \qquad \lim_{\varepsilon \downarrow 0} P(\tau \leq \varepsilon^{\alpha - 2 - \delta}) = 1.$$

PROOF. To show (3.29), we may assume that $\delta > 0$ is sufficiently small. Take $\theta \in (\tilde{\nu}, \frac{\alpha - 3}{2})$. Then, since $\alpha - 2 - \delta > 2\theta + 1$ (for sufficiently small $\delta$), we have $\lim_{\varepsilon \downarrow 0} P(\tau_5 \geq \varepsilon^{\alpha - 2 - \delta}) = 1$ from Lemma 3.7. However, on the event $B^{(\varepsilon)} := \{\tau_5 \geq \varepsilon^{\alpha - 2 - \delta}\}$, we see $\max_{1 \leq i \leq N} |x_i(t) - \bar{x}_i(t)| \leq \varepsilon^\theta$ at $t = \varepsilon^{\alpha - 2 - \delta}$, and therefore

$$\|\nabla \mathbf{h}(\mathbf{x}(t))\|_\infty = \max_{1 \leq i \leq N-1} |x_{i+1}(t) - x_i(t) - a|$$

$$\leq \max_{1 \leq i \leq N-1} |\bar{x}_{i+1}(t) - \bar{x}_i(t) - a| + 2 \max_{1 \leq i \leq N} |x_i(t) - \bar{x}_i(t)|$$

$$\leq (|b'_3 - a| \vee |b'_2 - a|) \{N \exp(-c_* \varepsilon^{-\alpha} N^{-2} \varepsilon^{\alpha - 2 - \delta})\}^{1/2} + 2\varepsilon^\theta \leq \varepsilon^{\tilde{\nu}}$$

if $\varepsilon > 0$ is sufficiently small. We have used (3.26) and then $\theta > \tilde{\nu}$ for the third line. This implies that $x(t) \in \mathcal{M}^{\nabla, N}(\varepsilon^{\tilde{\nu}})$ at $t = \varepsilon^{\alpha - 2 - \delta}$ and therefore $\tau \le \varepsilon^{\alpha - 2 - \delta}$ on the event $B^{(\varepsilon)}$, which proves (3.29). $\square$

We are now at the position to complete the proof of Theorem 3.1.

PROOF OF THEOREM 3.1. Combining Proposition 3.8 with Proposition 3.5 by means of the strong Markov property of $\mathbf{x}(t)$, we obtain

$$\lim_{\varepsilon \downarrow 0} P(\mathbf{x}(t) \in \mathcal{M}^{\nabla, N}(\varepsilon^{\tilde{\nu}}) \text{ for some } t \le \varepsilon^{-2-\delta} + \varepsilon^{\alpha - 2 - \delta}) = 1.$$

However, since $\varepsilon^{\alpha - 2 - \delta} \ll \varepsilon^{-2-\delta}$, the factor $\varepsilon^{\alpha - 2 - \delta}$ may be omitted by replacing $\delta$ if necessary. Hence, by introducing the time change (1.4), we obtain the conclusion. $\square$

We finally consider the case where the initial configuration $\mathbf{x}(0)$ consists of $n$ chains: $\mathbf{x}(0) = \mathbf{x}^{(1)}(0) \cup \cdots \cup \mathbf{x}^{(n)}(0)$ arranged from left to right with particles' numbers $N_1 \sim \rho_1 \varepsilon^{-1}, \ldots, N_n \sim \rho_n \varepsilon^{-1}$ and fluctuations $\varepsilon^{\nu}, \nu > 2 + (n-1)/2$, respectively, where $\rho_1, \ldots, \rho_n > 0$. Let $\mathbf{x}^{(\varepsilon)}(t) = (x_i^{(\varepsilon)}(t))_{i=1}^N, N = \sum_{\ell=1}^n N_\ell$ be the solution of the SDE (1.1) scaled macroscopically in time and starting at $\mathbf{x}(0)$. Denote the macroscopic center of mass of the associated $\ell$th rod by $\eta^{(\varepsilon, \ell)}(t) := \frac{\varepsilon}{N_\ell} \sum_{i=M_{\ell-1}+1}^{M_\ell} x_i^{(\varepsilon)}(t), 1 \le \ell \le n$, where $M_\ell = \sum_{\ell'=1}^{\ell} N_{\ell'}$ for $1 \le \ell \le n$ and $M_0 = 0$.

COROLLARY 3.9. Assume $\alpha > 2\nu + 3$. Then, the process $\{\eta^{(\varepsilon, \ell)}(t)\}_{\ell=1}^n$ converges to $\{\eta^{(\ell)}(t)\}_{\ell=1}^n$ as $\varepsilon \downarrow 0$ weakly in the space $C([0, T], \mathbb{R}^n)$ for every $T > 0$ if $\eta^{(\ell)}(0) = \lim_{\varepsilon \downarrow 0} \eta^{(\varepsilon, \ell)}(0)$ exist.

The limit process $\{\eta^{(\ell)}(t)\}_{\ell=1}^n$ of Corollary 3.9 is constructed as follows:

(1) $\{\tilde{\eta}^{(\ell)}(t) := \eta^{(\ell)}(t) - a(\sum_{\ell'=1}^{\ell-1} \rho_{\ell'} + \frac{1}{2} \rho_\ell)\}_{\ell=1}^n$ perform the Brownian motions with speeds inversely proportional to $\rho_\ell$ independently with each other until the time $\tau^{(1)} = \inf\{t; \tilde{\eta}^{(\ell)}(t) = \tilde{\eta}^{(\ell+1)}(t)$ for some $1 \le \ell \le n-1\}$.

(2) If the equality in the infimum for $\tau^{(1)}$ holds for $\ell = \ell^{(1)}$, then $\tilde{\eta}^{(\ell^{(1)})}(t) = \tilde{\eta}^{(\ell^{(1)}+1)}(t)$ for all $t \ge \tau^{(1)}$.

(3) The system $\{\tilde{\eta}^{(\ell)}(t); \ell \ne \ell^{(1)} + 1\}$ is afreshed at the time $\tau^{(1)}$ and, after $\tau^{(1)}$, each of them performs the Brownian motion with the same speed as above except $\ell = \ell^{(1)}$, for which the new speed is inversely proportional to $\rho_{\ell^{(1)}} + \rho_{\ell^{(1)}+1}$. The evolution continues independently until the time $\tau^{(2)} = \inf\{t; \tilde{\eta}^{(\ell)}(t) = \tilde{\eta}^{(\ell+1)}(t)$ for some $1 \le \ell(\ne \ell^{(1)}) \le n-1$ or $\tilde{\eta}^{(\ell^{(1)})}(t) = \tilde{\eta}^{(\ell^{(1)}+2)}(t)\}$.

(4) After the time $\tau^{(2)}$, the procedure is continued similarly along with the coagulation times $\tau^{(2)} < \tau^{(3)} < \cdots < \tau^{(n-1)}$. A single rod is finally left after the time $\tau^{(n-1)}$.

The proof of Corollary 3.9 is immediate from Theorems 2.2 and 3.1. Note that $\{\eta^{(\varepsilon,\ell)}(t)\}_{\ell=1}^{n}$ are independent Brownian motions until the time when the minimal distance between two of $n$ rods becomes $\varepsilon b$ and the coagulation of two rods occurs within the time interval of length $\varepsilon^{1-\delta}, \delta > 0$. The probability that more than three rods interact within the same such small time interval is negligible as $\varepsilon \downarrow 0$. We can therefore continue the argument given just after Theorem 3.1 also for $n \geq 2$.

## REFERENCES

[1] Funaki, T. (2004). Zero temperature limit for interacting Brownian particles. I. Motion of a single body. *Ann. Probab.* **32** 1201–1227. MR2060296

[2] Ikeda, N. and Watanabe, S. (1989). *Stochastic Differential Equations and Diffusion Processes*, 2nd ed. North-Holland, Amsterdam. MR1011252

[3] Itô, K. and McKean, H. P., Jr. (1974). *Diffusion Processes and Their Sample Paths*, 2nd printing, corrected. Springer, Berlin. MR345224

[4] Karatzas, I. and Shreve, S. E. (1991). *Brownian Motion and Stochastic Calculus*, 2nd ed. Springer, New York. MR1121940

[5] Lang, R. (1979). On the asymptotic behaviour of infinite gradient systems. *Comm. Math. Phys.* **65** 129–149. MR528187

[6] Lang, R. and Nguyen, X.-X. (1980). Smoluchowski's theory of coagulation in colloids holds rigorously in the Boltzmann–Grad-limit. *Z. Wahrsch. Verw. Gebiete* **54** 227–280. MR602510

[7] Mullins, W. W. (1992). A one dimensional stochastic model of coarsening. In *On the Evolution of Phase Boundaries* (M. E. Gurtin and G. B. McFadden, eds.) 101–105. Springer, New York.

Graduate School
of Mathematical Sciences
University of Tokyo
3-8-1 Komaba Meguro-ku
Tokyo 153-8914
Japan
e-mail: funaki@ms.u-tokyo.ac.jp